\documentclass[12pt]{article}
\usepackage{amssymb}
\usepackage{amsmath}
\usepackage{amsthm}
 \newcommand{\nn}{\mathbb N}
 \newcommand{\rr}{\mathbb R}
 \newcommand{\cc}{\mathbb C}

\newtheorem{thm}{Theorem}
\newtheorem{coro}{Corollary}

\newtheorem{prop}{Proposition}
\newtheorem{lem}{Lemma}

\begin{document}
\normalsize

\title{Complemented basic sequences in Frechet spaces with finite dimensional decomposition}
\author{ Hasan G\"ul\footnote{Middle East Technical University, e032209@metu.edu.tr} \ and S\"uleyman \"Onal\footnote{Middle East Technical University, 
osul@metu.edu.tr}}
\maketitle

\noindent
{\bf Abstract:} {\normalsize
Let $ E $ be a Frechet-Montel space and $ (E_{n})_{n \in \nn} $ be a finite dimensional unconditional decomposition of $ E $ with $ \dim(E_{n})\leq k $ for some fixed $ k \in \nn $ and for all $ n \in \nn $. Consider a sequence $ (x_{n})_{n \in \nn} \ $ formed by taking an element $ x_{n} $ from each $  E_{n} $ for all $ n \in \nn $. Then $ (x_{n})_{n \in \nn} $ has a subsequence which is complemented in $ E $.
}\\


\section{Introduction}
\paragraph{•}
Dubinsky has professed that approximation of an infinite dimensional space by its finite dimensional subspaces is one of the  basic tenets of Grothendieck's mathematical philosophy in creating nuclear Frechet spaces \cite{DUB2}. That approach encompasses bases, finite dimensional decompositions and the bounded approximation property of spaces and raises questions about conditions for which a space must satisfy in order to admit such structures and also about the relation of existence or absence of these structures with subspaces, quotient spaces,  complemented subspaces, etc. of the space.

He has ended his survey with a list of problems that were settled and open at that time \cite{DUB2}, that is, in 1989. For example, it has been shown that quotient spaces of a nuclear Frechet space may not have a basis \cite{DUM}, or similarly, that a nuclear Frechet space does not necessarily have a basis \cite{BES}. On the question of complemented subspace with basis, the following was the most general result known at that time \cite{DUB2}. A nuclear Frechet space with a decomposition into two dimensional blocks contains a complemented basic sequence \cite{DUB}, and a shorter proof of the same result by utilizing Ramsey Theorem was obtained by Ketonen  and Nyberg in \cite{KET} , which also includes a remark to the effect that the method of proof cannot be extended to spaces with decompositions into higher dimensional blocks. On the other hand, 
Krone and Walldorf showed that complemented subspaces with a strong finite-dimensional decomposition of nuclear K\"othe spaces have a basis \cite{KRO}.

Our first aim in this manuscript is to show that  nuclear Frechet spaces which have decompositions of $ k $ dimensional blocks contain a complemented basic sequence, where $ k \in \nn $. In fact, we have proved the following : Let $ E $ be a Frechet-Montel space and $ (E_{n})_{n \in \nn} $ be a finite dimensional decomposition of $ E $ with $ \dim(E_{n})\leq k $ for some fixed $ k \in \nn $ and for all $ n \in \nn $. Consider a sequence $ (x_{n})_{n \in \nn} \ $ formed by taking an element $ x_{n} $ from each $  E_{n} $ for all $ n \in \nn $. Then $ (x_{n})_{n \in \nn} $ has a subsequence which is complemented in $ E $.

After preliminaries in section 2, we establish the above result
 in section 3. 
\section{Preliminaries}
\paragraph{•}
 Recall that a locally convex (topological vector) space (over the field $\rr \ \text{\rm or} \
\cc) \ \ E $ is called a \textit{Frechet space} if it is a complete and metrizable.

If E is metrizable, then E has a countable basis of neighbourhoods of 0 and if E is also locally convex, we may choose a basis with elements that are absolutely convex, say $ \mathcal{B}= \lbrace U_{n} : n \in \nn \rbrace $ with $ U_{1} \supseteq U_{2} \supseteq \cdot\cdot\cdot $ . By setting $p_{n}( \cdot) : = \parallel \cdot \parallel_{U_{n}}$ we obtain an increasing sequence  $(p_{n})_{n \in \nn}$ of continuous
seminorms.  Thus, the topology of E can be given by an increasing sequence of continuous norms, a \textit{fundamental system of seminorms} .

A locally convex space  whose topology is defined using a countable set of compatible norms  i.e. norms such that if a sequence  that is fundamental in the norms  and  converges to zero in one of these norms, then it also converges to zero in the other, is called \textit{countably normed}. Note that for a countably normed space $ E $, $ (\hat{E}, \parallel \cdot\parallel_{k+1}) \rightarrow ( \hat{E}, \parallel \cdot\parallel_{k})$ , where $ \hat{E} $ denotes closure of $ E $, is an injection and all spaces we deal in this manuscript are countably normed.

  A subspace $ Y $ of a topological space $ X $ is called \textit{complemented} if $ Y $ is image of a continuous projection.

 A sequence $(x_{n})$ in a Frechet space $ E $ is called a \textit{basic sequence} if it is a Schauder basis of its closed linear span  $\overline{\text{ \rm span} \lbrace x_{n} : n \in \nn \rbrace} $   , which we also denote by $ [x_{n}]_{n \in \nn}$  .

A \textit{finite dimensional decomposition} (FDD in short) in $ E $ is a sequence $ (A_{n})_{n \in \nn} $ of continuous linear operators $  A_{n} : E \rightarrow E $ such that $ \dim(A_{n}(E)) < \infty $, $ A_{n}A_{m} = \delta_{mn}A_{n} $ and $x =\Sigma_{n=1}^{\infty} A_{n}x$ for all $x \in E $.

A $ k- $FDD is an FDD in which $ (A_{n})_{n \in \nn} $ can be chosen such that \\ $ \dim (A_{n}(E)) \leq k $ for some $ k \in \nn $. A \textit{strong FDD} is an FDD which is a k-FDD for some $ k \in \nn $.

Suppose that $ E $ has strong FDD property.

 Then we define $ E_{n} := A_{n}(E) $, for all $ n \in \nn $ and 
let $ M $ be an infinite subset of $ \nn $. Denote  $ E_{M} = \overline{\bigoplus_{n \in M} E_{n}} $ so that $ E= E_{\nn} $

 Let $ G $ be a subspace of $ E $. If $ G= \overline{\oplus_{n \in \nn}(G\cap E_{n})} $, then we call $ G $ a \textit{step subspace} of $ E $.

  We denote a sequence $ (x_{n})_{n \in \nn} \ $ formed by taking an element $ x_{n} $ from each $  E_{n} $ for all $ n \in \nn $ by $ (x_{E_{n}})_{n \in \nn} $. We define a quotient space of $ E $ obtained by 'dividing' $ E $ by a step subspace as \textit{  step quotient} of $ E $.

A complemented step subspace $ K $ of a Frechet space $ E $ is called \textit{naturally complemented} if complement of $ K $ can be made a step subspace. Equivalently, let $ P : E \rightarrow E$ be a projection with $ P(E) = K $. Then $ K $ is naturally complemented if $ P $ commutes with any projection map $ A_{n} $  defined above on $ E $.

$ \omega $ denotes the space of all sequences with real elements. We refer to \cite{Vogt} for any unexplained definitions and notations.
 
\section{Main result} 
The following lemmata are needed, and since the first one is well known, it is given without proof.

\begin{lem} \label{l1}
Let $ E $, $ F $ be topological vector spaces and $ T : E \rightarrow F $ be a continuous linear operator . Suppose restriction of $ T $ to a subspace $ M $ of $ E $ is an isomorphism and $ T(M) $ is complemented in $ F $. Then $ M $ is complemented in $ E $.

\end{lem}
 START\\
Although the following result, the existence of basic subsequence of certain sequences in a Banach spaces,
 is also known ( see, for example \cite{DST}), we include here with a sketch of a proof.
\begin{lem} \label{l2}
Let $ (x_{n}) $ be a normalized injective sequence  $ 
( x_{n} \neq x_{m}  $ for $ n \neq m )$ in a Banach 
space $ (E, \parallel \cdot \parallel) $ such that $ 
\underset{n \rightarrow \infty} \lim u(x_{n}) = 0 $ 
for all $ u \in B $ where $ B $ is a $ w^* $ dense , i.e. $\sigma(E, E^{'}) $ dense , subset of $ E^{'} $. Then, there exists an infinite subset $ M $ of $ \nn $ such that $ (x_{n})_{n \in M} $ is a basic sequence in $ E $.

\end{lem}

\begin{proof}

Since $ \parallel x_{n} \parallel = 1$ for all $ n \in \nn $ and $ x_{n} \rightarrow 0 $ in $\sigma(E, B) $, $ (x_{n}) $ does not have any $ \parallel \cdot \parallel $-Cauchy subsequence. Otherwise, if a subsequence converges to $x $ for some $ x \in E $, then  $ u(x)=0 $ for $ u \in B $ implies that $  \parallel x \parallel = 0  $ contradicting the fact that $ \parallel \cdot \parallel $ is continuous.

Thus, there exists an infinite subset $ M $ of $ \nn $ such that either $ (x_{n})_{n \in M} $ is equivalent to a basic sequence in $ l_{1} $ or $ (x_{n})_{n \in M} $ is weakly Cauchy with respect to  $\sigma(E, E^{'}) $ by Rosenthal's $ l_{1} $ theorem. If former is the case, then $ (x_{n}) $ has a basic subsequence.

For the latter case, we may assume, WLOG, $ (x_{n}) $ itself is a weakly Cauchy sequence. Now, suppose the bounded Cauchy sequence $ (x_{n}) $ converge to $ x \in E'' $ in $ \sigma(E'', E) $. If $ x \in E $, then $ x= 0 $ by the assumption that $ x_{n} \rightarrow 0 $ in $\sigma(E, B) $. Thus, $ (x_{n}) $ is a weakly null normalized sequence and by Bessaga-Pelczynski selection principle \cite{BP}, we can obtain a basic subsequence of $ (x_{n}) $. If $ x \in E'' \setminus E $, then the sequence $ (x-x_{n}) \rightarrow 0 $ in $ \sigma (E'', E') $ (see \cite{DST} Chapter IV, or Ex. 10 in page 55). Thus, we may obtain a basic subsequence of $ (x-x_{n})$. WLOG, let us assume that 
$ (x-x_{n})$ is basic. Recall that $ x \not\in [x_{n}] $ (here in this lemma, $ [\cdot] $ denotes the norm-closure of the span of $ \cdot $). Also, we claim that we can find $ k \in \nn $ such that 
$ x \not\in [x-x_{n}]_{n \geq k} $. (If $ x \in [x-x_{n}]_{n \geq 0} $ or $ x \in [x-x_{n}]_{n \geq 1} $, we set $ x = \underset{n \in \nn}\Sigma a_{n}(x - x_{n}) $ and if $ a_{n_{0}} $ is the first nonzero term, we let $ k = n_{0}+1 $.) Now let $ P_{1} $ and $ P_{2} $ be continuous projections on $[x-x_{n}]_{n \geq k}  $ and $ [x_{n}]_{n \geq k} $ with kernels (containing hg?) $ [x] $ and $ \kappa $ be a (the) basic constant (?) of  $[x-x_{n}]_{n \geq k}  $. We may assume $ k=1 $ by reindexing. Let $ m \in \nn $ such that$ \parallel \Sigma_{n=1}^{m}a_{n}x_{n}\parallel \leq 1 $. Then,  $ \parallel \Sigma_{n=1}^{m}a_{n}(x_{n}- x) + \Sigma_{n=1}^{m}a_{n}x \parallel \leq 1 $. \\ But $ \parallel \Sigma_{n=1}^{m}a_{n}(x_{n}- x)\parallel \leq \parallel P_{1} \parallel$ , the operator norm of $ P_{1} $. Thus, there exists $ l \in \nn $ such that $ \parallel \Sigma_{n=1}^{l}a_{n}(x_{n}- x)\parallel \leq \parallel P_{1} \parallel \cdot\kappa$ (?). Hence, $ \parallel \Sigma_{n=1}^{l}a_{n}x_{n}-\Sigma_{n=1}^{l} x\parallel \leq \parallel P_{1} \parallel \cdot\kappa$ . Therefore,  $ \parallel \Sigma_{n=1}^{l}a_{n}x_{n}\parallel \leq \parallel P_{1} \parallel \cdot\kappa \cdot \parallel P_{2} \parallel$ which shows that $ (x_{n})_{n \in \nn} $ is a basic sequence with basic constant $\parallel P_{1} \parallel \cdot\kappa \cdot \parallel P_{2} \parallel$ (?). Since $ x_{n}\in E $, $ [x_{n}] \subseteq E$, which is to say, $ x_{n}$ is a basic sequence in $ E $.

\end{proof}

In the above proof, we refer to Rosenthal $ l_{1} $-
theorem, but actually we do not need to do that. 
Since $ x_{n} \succ (?)x $ in $ \sigma (E'', E') $, $ 
y \in \underset{k \geq n} \bigcap [x_{m}]_{m \geq k}$ 
with respect to $ \sigma (E'', E') $. Hence, $ [x_{n}-x] \succ 0 $ with respect to $ \sigma (E'', E') $. Thus a similar approach as we claimed $ \sigma (E, B) $ -converging to $ 0 $ by ....

We prefer that approach since in a countably normed Frechet space $ E $, a subset $ B $ with the above property (ies) can easily be found. \\(END)\\

It has been shown that the restriction of unbounded linear maps to infinite dimensional subspaces between some classes of Frechet spaces is an isomorphism. \cite{TY}

We shall show that for any sequence $ (x_{n}) $ in a Frechet space $ E $ and any linear map $ T $, $ (T(x_{n})) $ either has a subsequence, when normalized, is bounded or restriction of $ T $ on $ [x_{n}]_{n \in \nn}$ is an isomorphism.

\begin{prop}

Let $ E $ and $ F $ be Frechet spaces, $ T : E \rightarrow F $ be a linear map and $ (x_{n}) $ be a sequence in $ E $. Suppose $ E $ has a continuous norm and $ F $ is countably normed. ( spaces with basis, their subspaces etc. are c. n. ) Then there exists an infinite subset $ M $ of $ \nn $ such that

(i) $ ( \frac{1}{\parallel x_{n} \parallel} Tx_{n} )_{n \in M}$ is a bounded subset of $ F $ for some continuous norm $ \parallel\cdot\parallel $ on $ E $,

or

(ii) restriction of  $T$  to $ [x_{n}]_{ n \in M}  $, that is,   $ T\mid_{ [x_{n}]_{n \in M}} $ is an isomoporhism and  furthermore $ [x_{n}]_{ n \in M} =\lambda(A) $ for some nuclear K\"othe space $ \lambda(A) $. 

\end{prop}

\begin{proof}

Let $ ( \parallel \cdot \parallel_{n} ) $ and  $ ( \mid \cdot \mid_{n} ) $ be increasing sequences of norms which define topologies of $ E $ and $ F $, respectively, with $ \mid Tx \mid_{k} \leq \parallel x \parallel_{k}$ and  $ F $ is countably normed. Without loss of generality, we may assume that $ T_{\mid [x_{n}]_{n \in \nn}} $ is an injection.

 Set $ M_{0}=\nn $. Consider the sequence $ ( \frac{Tx_{n}}{\parallel x_{n} \parallel_{1}} )_{n \in M}$. If it is bounded, we may take $ M $ to be $ \nn $ and $ k =1 $ to claim (i) holds.

Otherwise, $ ( \frac{Tx_{n}}{\parallel x_{n} \parallel_{1}} )_{n \in M_{0}}$ is unbounded

 Since$ ( \hat{F}, \mid \cdot \mid_{k+1}) \rightarrow ( \hat{F}, \mid \cdot \mid_{k}) $ is an injection for all $ k \in \nn $ there exists $ l_{1} \in \nn $ such that $ \sup (\frac{\mid Tx_{n}\mid_{l_{1}}}{\parallel x_{n}\parallel_{1}}) = \infty $.

Thus, we can find an infinite subset $ M_{1} $ of $ M $ such that $ \underset{n \in M_{1}}\Sigma \frac{\parallel x_{n}\parallel_{1}}{\mid Tx_{n}\mid_{l_{1}}} < \infty $ Since $ \mid Tx_{n}\mid_{1} \leq \parallel x_{n} \parallel_{1} $, the sequence $ (\frac{ Tx_{n}}{\mid Tx_{n}\mid_{l_{1}}})_{n \in M_{1}}  $ converges to zero. Thus, the normalized sequence $\frac{Tx_{n}}{\mid Tx_{n}\mid_{l_{1}}}$ converges to zero with respect to the topology $ \sigma(F, F_{1}^{'}) $ where $F_{1}^{'}  $ denotes the dual of the normed space $ ( F, \mid \cdot \mid_{1} )$. But  $F_{1}^{'}  $ is $ \sigma(F_{l_{1}}^{'}, (\hat{F}_{l_{1}}, \mid \cdot \mid_{l_{1}})) $ dense. By applying Lemma 2, that is, considering $ E $ to be $ (\hat{F}_{l_{1}}, \mid \cdot \mid_{l_{1}})) $ and the sequence to be  $ (\frac{ Tx_{n}}{\mid Tx_{n}\mid_{l_{1}}})_{n \in M_{1}}  $ following the notation of Lemma 2, we obtain a basic subsequence of $ (x_{n})_{n \in M} $.

Now consider the sequence $ ( \frac{Tx_{n}}{\parallel x_{n} \parallel_{l_{1}}} )_{n \in M_{1}}$. If it is bounded, by taking  $ M $ to be $ M_{1} $ and $ k =l_{1} $ we may  claim (i) holds.

Otherwise, there exists $ l_{2} \in \nn $ such that $ \lim_{n \rightarrow \infty} \frac{\mid Tx_{n} \mid_{l_{2}}}{\parallel x_{n}\parallel_{l_{1}}} = \infty $.  Similar to above, we may find an infinite subset $ M_{2} $ of $ M_{1} $ such that  the sequence $ (\frac{ Tx_{n}}{\mid Tx_{n}\mid_{l_{2}}})_{n \in M_{2}}  $ converges to zero and we may use Bessaga-Pelczynski selection principle \cite{BP} to obtain an basic sequence in $ (F, \parallel\cdot\parallel_{l_{2}}) $ if necessary. Note that we may choose $ \mid Tx\mid_{k}\leq\parallel x\parallel_{k} $  for all $ k \in \nn $, so that $ l_{2}>l_{1} $.

If we proceed in this manner, use Bessaga-Pelczynski selection principle \cite{BP} to obtain an basic sequence in $ (F, \parallel\cdot\parallel_{l_{k}}) $, if necessary, and assuming (i) does not hold in each step, we obtain an infinite decreasing sequence of $ ( M_{t})  $ of infinite subsets of $ \nn $. By applying diagonalization argument to $ ( M_{t})  $, we get an infinite subset $ M $ of $ \nn $ such that $ M \setminus M_{t} $ is finite for all $ t \in \nn $. Now, it is easy to see that the restriction of $ T $ to $ [x_{n}]_{n \in M} $ is an isomorphism and $ [x_{n}]_{n \in M} $ is nuclear that is, (ii) holds.

\end{proof}

\begin{coro} \label{c1}
   Let $ E $ be a metrizable locally convex space with a basis $  (x_{n})_{n \in \nn} $ with increasing continuous norms $ \parallel\cdot\parallel_{k} $ for all $ k \in \nn $ and $ F $ be a countably normed Frechet space with increasing continuous seminorms $ \mid\cdot\mid_{t} $ for all $ t\in \nn $ which gives topologies of the respective spaces. Let the identity operator $ I_{F} $ on F,   have a decomposition $ I_{F} = S_{1}+S_{2}+\cdot\cdot\cdot S_{l}$ into some continuous linear maps $ S_{i} : F \rightarrow F $ for all $ 1 \leq i \leq l $ for some $ l \in \nn. $ Let $T : E \rightarrow F $ be a continuous linear operator from $ E $ into $ F $. Then there exists an infinite subset $ M $ of $ \nn $ such that either $$ (a) \  \left( \frac{Tx_{n}}{\parallel x_{n}\parallel_{k}}\right) _{n \in M}  \textit{   is a bounded subset of F for some \ } k \in \nn,$$  or \\
  
 (b) restriction of  $S_{i}T$  to $ [x_{n}]_{ n \in M}  $ is an isomorphism for some  $i \in \lbrace 1,2,...,l \rbrace $.

\end{coro}

\begin{proof}  WLOG, we may assume $ I_{F} = S_{1} + S_{2} $. Apply lemma \ref{l2} to $ S_{1}T $ and suppose the conclusion (b) does not hold. Then, $ ( \frac{S_{1}Tx_{n}}{\parallel x_{n} \parallel_{1}} )_{n \in M}$ is bounded for some infinite subset $ M $ of $ \nn $. Now apply  lemma \ref{l2} to the restriction of $ S_{2}T$ to $[(x_{n})]_{n \in M} $. If (i) in lemma \ref{l2} holds, then conclusion (a) holds.
Otherwise (b) holds.

\end{proof}

\begin{prop} \label{p1} Let $ F $ be a Frechet space with an unconditional basis $ (e_{i})_{i \in \nn} $, satisfying $ e^{\ast}_{i}(e_{j})= \delta_{ij} $ for $ e_{i}^{\ast} \in F' $, which admits a continuous norm. Let $ (x_{n})_{n \in \nn} $ be a sequence of nonzero vectors in $ F $ such that the number of elements of the set $ A_{n} = \lbrace i \in \nn : e_{i}^{\ast}(x_{n}) \not = 0 \rbrace $ is at most $ k $ for each $ n \in \nn $ for some $ k \in \nn $ and $ A_{n}\cap A_{m} = \emptyset $ when $ n \not= m $ . If $ F $ is Montel or isomorphic to some K\"othe space $ \lambda^{p}(A) $ with the canonical base mentioned above, then there is a complemented subsequence $( x_{n_{m}})_{m \in \nn} $ of  $ (x_{n})_{n \in \nn} $.

\end{prop}

\begin{proof} Let $ ( L_{i})_{i=1}^{k} $ be a partition of $ \nn $ such that $ L_{i}\cap A_{n} $ has at most one element for each $ n \in \nn $ and $ 1 \leq i \leq k $. Set $ F_{i} := [e_{n}]_{n \in L_{i}} $ and $ E := [x_{n}]_{n \in \nn} $ and $ T : E \rightarrow F := \oplus_{i=1}^{k}F_{i} $ be the inclusion map. Then by Lemma 2, there exists an infinite subset $ M $ of $ \nn $ such that either (i) $  \  \left( \frac{x_{n}}{\|( x_{n})\|}\right) _{n \in M}$  is bounded where $ \| \cdot \| $ is a continuous norm on $ F $ or (ii) the restriction of natural projection $ P_{j} $ from $ F $ onto $ F_{j} $ is an isomorphism of $ [x_{n}]_{n \in M} $ for some $ j \in \lbrace 1,2,...,k \rbrace $.

If (ii) holds then $ [x_{n}]_{n \in M} $ is complemented in $ F $ by Lemma 1. 

Now, if $ F $ is Montel, then (i) does not hold. Therefore, $ (x_{n})_{n \in M} $ is a complemented subsequence  of  $ (x_{n})_{n \in \nn} $.

If $ F \simeq \lambda^{p}(A) $ and (i) holds then $ [x_{n}]_{n \in \nn} $ is isomorphic to $ l_{p} $ and
$ P(x) := \Sigma_{n \in M}^{}x_{n}^{\ast}(x)x_{n} $ is a projection from $ F $ onto  $ [x_{n}]_{n \in M} $ where $ x^{\ast}_{n}(x_{m})= \delta_{nm} $ for $ x_{n}^{\ast} \in F' $ and $\parallel x_{n}^{\ast} \parallel \leq r $ for each $ n,m \in M $ and for some $ r \in \rr $.
\end{proof}

\begin{coro} \label{c2}
Let $ (E_{n})_{n \in \nn} $ be a two dimensional decomposition of a nuclear Frechet space E. Then each sequence of nonzero vectors $ (x_{E_{n}})_{n \in \nn} $ has a complemented subsequence.

\end{coro}

\begin{proof}   Since $ E $ is a nuclear Frechet space with 2-FDD, there exists a sequence $ (y_{E_{n}})_{n \in M} $ for some infinite subset $ M $ of $ \nn $ such that $[y_{E_{n}}]_{n \in M} $ has a complementary basis  \cite{DUB}, which can be made naturally complemented by , say   $[z_{E_{n}}]_{n \in M} $.  Thus, we may write $ E_{n}= [y_{E_{n}}]+[z_{E_{n}}] $ and $ E_{M}=   [y_{n_{i}}]_{i \in M} \oplus  [z_{n_{i}}]_{i \in M} $ .

 Let $ P $ be the projection map on $ E $ with range $ [y_{E_{n}}]_{n \in M} $. Then $$ P^{-1}(0)= [z_{E_{n}}]_{n \in M} \bigoplus \oplus_{n \in \nn \setminus M} E_{n}.$$ Recall  that $ E_{M} := [E_{n}]_{n \in M} $ and note $ I_{M}= P\mid _{M} + (I - P)\mid_{M}$.  By applying proposition \ref{p1} to the space $ [x_{n}]_{_{n \in M}} $and taking $ I_{M} $ to be $ T $,  we obtain the fact that 
$ x_{E_{n}} $ has a complemented subsequence.

\end{proof}

Let $ E $ be a locally convex space such that $ E = 
\overline{\oplus_{n \in \nn} E_{n}} $ where $ \dim(E_{n}) = k $ 
for each $ n \in \nn $ and for some $ k \in \nn $ We 
say that $ E $ satisfies \textit{condition (*)}, if $ F 
$ is 
any 2-dimensional step quotient of step subspace  of $ E $ with decomposition $ F = \oplus_{n \in \nn} 
F_{n}$ contains closed subspaces $ G_{x} $ for any $ (x_{F_{n}})_{n \in \nn} $ such that
$F_{M} = [x_{F_{n}}]_{n \in M}+ G_{x} $ and $ \dim (G_{x} \cap F_{n})=1 $ for some infinite subset $ M $ of $ \nn  $ and for all $ n \in \nn $.

\begin{prop} \label{p2}    Let $ E$ be a locally convex space such that $ E = \overline{\oplus_{n \in \nn}E_{n}} $ where $ \dim(E_{n}) = k $ for each $ n \in \nn $, for some $ k \in \nn $, and satisfy condition (*). Then, there are closed subspaces $ G_{1},G_{2},\cdot\cdot\cdot,G_{k} $ of $ E $ and an infinite subset $ M $ of $ \nn $. such that 
$$ E_{M} = G_{1}+ G_{2}+ \cdot\cdot\cdot +G_{k} \text{ \rm \ \ where } \dim(G_{i}\cap E_{n})= 1        \text{                                           \rm for all } n \in M , \text{                                           \rm  } 1\leq i \leq k .$$ 

\end{prop}

\begin{proof} Let $ E$ be a locally convex space such that $ E = \overline{\oplus_{n \in \nn}E_{n}} $ where $ \dim(E_{n}) = k $ for each $ n \in \nn $, for some $ k \in \nn $, and satisfy condition (*). We proceed with induction on $ k $. For the case $ k=1 $, it is immediate.

 Let $ (x_{E_{n}})_{n \in \nn} $ be nonzero vectors. Let $ Q : E \rightarrow E / [x_{E_{n}}] =: F $ be the quotient map. Thus, $ F $ satisfies condition (*) with $ \dim(F_{n}) = k - 1 $. By induction hypothesis, we have closed subspaces $ G_{1},G_{2},\cdot\cdot\cdot,G_{k-1} $ of $ F $ such that $$ F_{M} = G_{1}+ G_{2}+ \cdot\cdot\cdot +G_{k-1} \text{ \rm \ \ where } \dim(G_{i}\cap F_{n})= 1        \text{                                           \rm for all } n \in \nn , \text{                                           \rm  } 1\leq i \leq k-1 ,$$ and for some infinite subset $ M $ of $ \nn $. Set $ H_{n}^{i}:= Q^{-1}(G_{i})\cap E_{n} $ which are 2-dimensional subspaces of $ E $ containing $ x_{E_{n}}  $ for  all $ n \in M $ and for $ 1 \leq i \leq k-1 $.

Define $ \mathcal{H}_{i} := [H_{n}^{i}]_{n \in M} $. Then  $ \mathcal{H}_{i}$ is 2-FDD step subspace of $ [E_{n}]_{n \in M} $, a space which satisfy condition (*) for $ 1 \leq i \leq k-1 $. If we consider $ (x_{E_{n}})_{n \in M} $ and $\mathcal{H}_{1}  $ in conjunction with condition (*), we conclude that there exists a closed subspace $ \mathcal{G}_{1} $ of $ E_{M} $ such that $ [\mathcal{H}_{1}]_{M_{1}} = [x_{E_{n}}]_{n \in M_{1}} +\mathcal{G}_{1}  $ for some infinite subset $ M_{1} $ of $ M $.

 Now consider $ (x_{E_{n}})_{n \in M} $ and $\mathcal{H}_{2}  $ to obtain a closed subspace $ \mathcal{G}_{2} $ of $ \mathcal{H}_{1}$ such that
$ [\mathcal{H}_{2}]_{M_{2}} = [x_{E_{n}}]_{n \in M_{2}} +\mathcal{G}_{2}  $ for some infinite subset $ M_{2} $ of $ M_{1} $.

After proceeding in this manner, we end up with 
$$ E_{L}= [x_{E_{n}}]_{n \in L}\oplus \mathcal{G}_{1} \oplus \mathcal{G}_{2} \oplus \cdot\cdot\cdot \oplus \mathcal{G}_{k-1} ,$$ where $ L:= M_{k-1}\subseteq M_{k-2}\subseteq \cdot\cdot\cdot \subseteq M_{1}\subseteq M \subseteq \nn.$

\end{proof}

\begin{prop} \label{p3}
Let $ E $ be a Frechet space with a k-dimensional unconditional decomposition $ (E_{n})_{n \in \nn} $ and suppose $ E $ satisfies condition (*). Then any sequence $ (x_{E_{n}})_{n \in \nn} $ has a naturally complemented subsequence in $ E $.

\end{prop}

\begin{proof} Let $ E $ be a Frechet space with a k-dimensional unconditional decomposition $ (E_{n})_{n \in \nn} $ and suppose $ E $ satisfies condition (*).   Then, there are closed subspaces $ G_{1},G_{2},\cdot\cdot\cdot,G_{k} $ of $ E $ and an infinite subset $ M $ of $ \nn $ such that 
$$ E_{M} = G_{1}+ G_{2}+ \cdot\cdot\cdot +G_{k} \text{ \rm \ \ where } \dim(G_{i}\cap E_{n})= 1        \text{                                           \rm for all } n \in M) , \text{                                           \rm  } 1\leq i \leq k ,$$  and $ G_{1}= [x_{n}]_{n \in M} $ by proposition \ref{p2}. Graph of projections $ P_{i} $ with ranges $ G_{i} $, respectively,  are closed, thus, by closed graph theorem each projection $ P_{i} $ is continuous, for $ 1 \leq i \leq k $. Thus, by lemma \ref{l1}, $ E_{M} $ is complemented.

Therefore, any sequence of nonzero vectors $ (x_{E_{n}})  $ has a complemented subsequence.

\end{proof}

The next result follows from Corollary 2 and Proposition 3.

\begin{thm} \label{t1}
 Let $ E $ be a nuclear Frechet space  with a k-dimensional decomposition $( E_{n})_{n \in \nn} $. Then any sequence of nonzero vectors $ x_{n} \in E_{n} $ has a complemented subsequence and the subspace $ [E_{n}]_{n \in M} $ of $ E $ has a basis for some infinite subset $ M $ of $ \nn $.

\end{thm}

\section{Some generalizations}

So far, we have dealt with nuclear Frechet spaces. Now we show that the result obtained in Theorem 1 can be generalized to Montel-Frechet spaces and then we turn our attention to Frechet spaces with unconditional FDD which does not admit a continuous norm.

\begin{prop} \label{p3} Let $ E =  [E_{n}]_{n \in \nn}  $ be a Montel-Frechet space with unconditional k-FDD with continuous norm. Then any sequence  $ (x_{E_{n}}) $ has a naturally complemented subsequence. Furthermore, there exists an infinite subset $ M $ of $ \nn $ such that $ E_{M} = \lambda(A) $  for some nuclear Kothe space $ \lambda(A) $.

\end{prop}

\begin{proof} It is sufficient to show that $ E_{M} $ is nuclear for some $ M \subseteq \nn $ by Theorem 1.

Let $ x_{E_{n}} $ be a sequence of nonzero vectors. By applying Lemma 2 to $ I_{E} $ , we end up the second conclusion in the lemma 2 holds, that is, $[x_{E_{n}}]_{L \in \nn}   $ is nuclear for some infinite subset $ L $ of $ \nn $, for $ E $ is Montel. We proceed by induction on $ k $.

For $ k=1 $, it is immediate.

Now define $ [{F_{n}}]_{n \in L} := [{E_{n}}]_{n \in L} / [x_{E_{n}}]_{L \in \nn}$. $ [{F_{n}}]_{n \in L}  $ is a Montel space, since $[{E_{n}}]_{n \in L}$ is Montel and $[x_{E_{n}}]_{L \in \nn}$ is nuclear. Also, note that $ [{F_{n}}]_{n \in L}  $ has $ k-1-$FDD and is with continuous norm. By induction hypothesis, there exists an infinite subset $ M $ of $ \nn $ such that $ F_{M}  $  is  nuclear.

Therefore, $ E_{M}  $ is nuclear by the three space property of nuclearity.

\end{proof}

   Now suppose $ E $ is not Montel but has an absolute p-basis instead of an unconditional one with k-FDD and continuous norm. Then, one can show that any sequence  $ (x_{E_{n}}) $ has a naturally complemented subsequence in that space. In fact, there is M, such that  either $ E_{M}=\lambda(A) $ or $ E_{M}=l_{p} $ or $ E_{M}=\lambda(A) \times l_{p}$ for $ 1 \leq p < \infty $ or $ p=0 $, that is, $ c_{0} $ in view of Lemma 2.

   If  $ E $ is a subspace of $ l_{p}^{\nn} $ or $ c_{0}^{\nn} $ with strong FDD , using a modification of  Bessaga-Pelczynski selection principle for the finite rank operators in Frechet spaces, we may suppose $ E $ to have an absolute p-basis and proceed as in the above paragraph.

   As for spaces which do not admit continuous norms,  we have the following which can be considered as a similar result for spaces with FDD instead of spaces with unconditional basis which was obtained by Floret and Moscatelli in \cite{FLO}. 

\begin{prop} \label{p5} Let $ E $ be a Frechet space with FDD and with a sequence of projections $ \lbrace P_{n} : n \in \nn \rbrace $ and a fundamental system of seminorms $ \parallel\cdot\parallel_{k} $. WLOG, we may assume $ \parallel P_{n}(x) \parallel_{k} \leq \parallel x \parallel_{k}$ for all $ k \in \nn $
and for all $ x \in E $. Then, either $ E $ has a continuous norm or $ E =\oplus E_{M_{i}} $ unconditionally where each $ E_{M_{i}} $ admits continuous norms for $ M_{i} \subseteq \nn  $ for all $ i \in \nn. $ 
\end{prop}

\begin{proof} Let $ E_{n} := P_{n}(E) $ for all $ n \in \nn $. Define 
$$ \varphi(n):=  min \lbrace k : \ \parallel \cdot \parallel_{k} \text{ \rm is a norm on} \ E_{n} \rbrace  .$$ Then, set $ M_{k} := \lbrace n \in \nn : \varphi(n)=k \rbrace . $ Now, either $ \nn = \bigcup_{k=1}^{L}M_{k} $ for some $ L \in \nn $ or $ \nn = \bigcup_{k=1}^{\infty}M_{k} $. In the former case, $ E $ admits a continuous norm. In the latter case, $E = \oplus E_{M_{k}}$  with $ \parallel \cdot \parallel_{k} $ is a norm on $E_{M_{k}}$. Note that if
 each $ M_{k} $ is finite then $ E = \omega $.

\end{proof}

If $ E \not=\omega $ but have a FDD with unconditional basis, then there exists an infinite subset $ M $ of $ \nn $ such that $ E_{M} $ admits a continuous norm. Furthermore, if $ E $ has a strong FDD, the above results are applicable for $ E_{M} $.

Hasan G\"ul,                 \ \hspace{50mm}             S\"uleyman \"Onal

Department of Mathematics,   \ \hspace{18mm} Department of Mathematics,

METU,                              \ \hspace{58mm}METU,

In\"on\"u Bulvari 06531          \ \hspace{36mm}In\"on\"u Bulvari 06531

Ankara, Turkey                  \ \hspace{44mm}Ankara, Turkey   

\begin{thebibliography}{HD}



\bibitem[1]{BES}
Bessaga C. and Dubinsky E.,
Nuclear Frechet spaces without bases. \\
Arch. Math. 31, 597-604. 6 (1978).


\bibitem[2]{BP}
Bessaga C. and Pelczynski A.,
On bases and unconditional convergence of series in Banach spaces.
Studia Math. 17, 151-164. (1958).


\bibitem[3]{DST}
Diestel J.,
Sequences and series in Banach spaces. \\
Springer-Verlag, New York. (1984) .



\bibitem[4]{DUB2}
Dubinsky E.,
Approximation properties of nuclear Frechet spaces,  
in Terzio\^glu T. (ed.),
Advances in the Theory of Frechet Spaces. Kluwer Academic Publishers, Dordrecht, 1-10. (1989) .




\bibitem[5]{DUB}
Dubinsky E.,
Complemented basic sequences in nuclear Frechet spaces with finite dimensional decomposition.
Arch. Math. 38, 138-150. (1982).


\bibitem[6]{DUM}
Dubinsky E. and Mityagin B. S.,
Quotient spaces without basis in nuclear Frechet spaces.
Canad. J. Math. 30, 1296-1305. (1978).

\bibitem[7]{FLO}
Floret K. and Moscatelli V. B.,
Unconditional basis in Frechet spaces.\\
Arch. Math.  47, 129-130. 2(1986).


\bibitem[8]{KET}
Ketonen T.  and Nyberg K.,
The Ramsey theorem and complemented basic sequences in Frechet spaces with absolute  finite dimensional decomposition.
Arch. Math. 53, 399-404. (1989).


\bibitem[9]{KRO}
Krone J.  and Walldorf V.,
Complemented subspaces with a strong finite-dimensional decomposition of nuclear K\"othe spaces have a basis. \\
Studia Math. 127, 1-7. 1(1998).

\bibitem[10]{TY}
Terzioglu T.  and Yurdakul M.,
Restrictions of unbounded continuous linear operators  on Frechet spaces\\
Arch. Math. 46, 547-550. (1986).

\bibitem[11]{Vogt}
Vogt D. and Meise R.,
Introduction to Functional Analysis.\\
Oxford University Press, New York. (2004).





\end{thebibliography}
\end{document}